\documentclass[a4paper,12pt,final]{amsart}
\usepackage{times,a4wide,mathrsfs,amssymb,dsfont, enumerate}

\newcommand{\C}{\mathbb{C}}
\newcommand{\Z}{\mathbb{Z}}

\newcommand{\QQ}{\mathbb{Q}}
\newcommand{\NN}{\mathbb{N}}
\newcommand{\PP}{\mathbb{P}}

\newcommand{\OO}{\mathcal O}
\newcommand{\Ss}{\mathcal S}

\newcommand{\XX}{\mathcal X}
\newcommand{\CCC}{\mathcal C}

\newcommand{\NNN}{\mathcal N}

\newcommand{\JJ}{\mathcal J}

\newcommand{\MM}{\mathcal M}
\newcommand{\FF}{\mathcal F}
\newcommand{\EE}{\mathcal E}
\newcommand{\PPP}{\mathcal P}
\newcommand{\AAA}{\mathcal A}

\newcommand{\gr}{\hbox{Gr}}
\newcommand{\wt}{\widetilde}
\newcommand{\ima}{\hbox{Im}}

\newcommand{\ide}{\hbox{id}}

\newtheorem{theorem}{Theorem}[section]
\newtheorem{claim}[theorem]{Claim}
\newtheorem{lemma}[theorem]{Lemma}

\newtheorem{proposition}[theorem]{Proposition}
\newtheorem{conjecture}[theorem]{Conjecture}
\newtheorem{nonumbering}{Theorem}

\newtheorem{convention}{Conventions}

\theoremstyle{definition}
\newtheorem{remark}[theorem]{Remark}
\newtheorem{definition}[theorem]{Definition}

\newtheorem{notation}[theorem]{Notation}

\newtheorem{nonumberingt}{Acknowledgements}

\begin{document}
\author[Robert Laterveer]
{Robert Laterveer}

\address{Institut de Recherche Math\'ematique Avanc\'ee,
CNRS -- Universit\'e 
de Strasbourg,\
7 Rue Ren\'e Des\-car\-tes, 67084 Strasbourg CEDEX,
FRANCE.}
\email{robert.laterveer@math.unistra.fr}

\title{On the Chow ring of some Lagrangian fibrations}

\begin{abstract} Let $X$ be a hyperk\"ahler variety admitting a Lagrangian fibration. Beauville's ``splitting property'' conjecture predicts that
fibres of the Lagrangian fibration should have a particular behaviour in the Chow ring of $X$. We study this conjectural behaviour for two very classical examples of Lagrangian fibrations.
\end{abstract}



\maketitle

\section{Introduction}

  For a smooth projective variety $X$ over $\C$, let $A^i(X):=CH^i(X)_{\QQ}$ denote the Chow groups (i.e. the groups of codimension $i$ algebraic cycles on $X$ with $\QQ$-coefficients, modulo rational equivalence). 
%Let $A^i_{hom}(X)$ denote the subgroup of homologically trivial cycles. 

The domain of algebraic cycles is an alluring treasure trove for anyone looking for open problems \cite{B}, \cite{J2}, \cite{Kim}, \cite{Mur}, \cite{Vo}, \cite{MNP}. Inside this treasure trove, one niche of particular interest is occupied by hyperk\"ahler varieties (i.e. projective irreducible holomorphic symplectic manifolds \cite{Beau1}, \cite{Beau0}). For these varieties, recent years have seen an intense amount of new constructions and significant progress in the understanding of their Chow groups \cite{Beau3}, \cite{V13}, \cite{V14}, \cite{V17}, \cite{SV}, \cite{V6}, \cite{Rie}, \cite{Rie2}, \cite{LFu}, 
\cite{Lin}, \cite{Lin2}, \cite{SYZ}, \cite{FTV}, \cite{FT}, \cite{MN}, \cite{V20}.
%Here, motivated by results for $K3$ surfaces and for abelian varieties, in recent years significant progress has been made in the understanding of Chow groups. 
%\cite{BV}, \cite{V13}, \cite{V14}, \cite{V17}, \cite{Fer}, \cite{SV}, \cite{V6}, \cite{Rie}, \cite{Rie2}, \cite{LFu2}, \cite{Lin}, \cite{Lin2}, \cite{FTV}.
Much of this progress has centered around the following conjecture:

\begin{conjecture}[Beauville, Voisin \cite{Beau3}, \cite{V17}]\label{conjbv} Let $X$ be a hyperk\"ahler variety. Let $D^\ast(X)\subset A^\ast(X)$ denote the $\QQ$-subalgebra generated by divisors and Chern classes of $X$. Then the cycle class maps induce injections
  \[ D^i(X)\ \hookrightarrow\ H^{2i}(X,\QQ)\ \ \ \forall i\ .\] 
  \end{conjecture}
  
  (For some cases where Conjecture \ref{conjbv} is satisfied, cf. \cite{Beau3}, \cite{V17}, \cite{BV},  \cite{Rie2}, \cite{LFu}, \cite{V6}, \cite{Yin}, \cite{FT}, \cite{LV2}, \cite{FLV2}.)
  
  The ``motivation'' underlying Conjecture \ref{conjbv} is that for a hyperk\"ahler variety $X$, the Chow ring $A^\ast(X)$ is expected to have a bigrading $A^\ast_{[\ast]}(X)$, where the piece $A^i_{[j]}(X)$ corresponds to the graded $\gr^j_F A^i(X)$ for the conjectural Bloch--Beilinson filtration. In particular, it is expected that the subring $A^\ast_{[0]}(X)$ injects into cohomology, and that $D^\ast(X)\subset A^\ast_{[0]}(X)$.
  
 In addition to divisors and Chern classes, what other cycles should be in the subring $A^\ast_{[0]}(X)$ (assuming this subring exists) ?
 A conjecture of Voisin provides more candidate members:
 
 \begin{conjecture}[Voisin \cite{V14}]\label{conjv} Let $X$ be a hyperk\"ahler variety of dimension $n=2m$. Let $Z\subset X$ be a codimension $i$ subvariety swept out by $i$-dimensional constant cycle subvarieties. There exists a subring $A^\ast_{[0]}(X)\subset A^\ast(X)$ injecting into cohomology, containing $D^\ast(X)$ and 
   \[ Z\ \ \in A^i_{[0]}(X)\ .\]
   \end{conjecture}
   
    A {\em constant cycle subvariety\/} is a closed subvariety $T \subset X$ such that the image of the natural map $A_0(T)\to A^n(X)$ has dimension $1$. In particular, Conjecture \ref{conjv} stipulates that Lagrangian constant cycle subvarieties (i.e., constant cycle subvarieties of dimension $m$) should lie in $A^m_{[0]}(X)$.
Some results towards Conjecture \ref{conjv} can be found in \cite{V13}, \cite{Lin2}, \cite{FLV}, \cite{FLV2}.

Amongst hyperk\"ahler varieties, of particular interest are those admitting a {\em Lagrangian fibration\/} (i.e. a proper surjective morphism $\pi\colon X\to B$ with connected fibers and $0<\dim B<\dim X$; in this case the general fiber of $\pi$ is an abelian variety that is Lagrangian with respect to the symplectic form on $X$ \cite{Mat}. In dimension 2, a Lagrangian fibration is an elliptic K3 surface\footnote{For background on Lagrangian fibrations, cf. the foundational \cite{Mat} as well as the recent \cite{BK}, \cite{HX} and the references given there.}).
As explained in \cite{mmj}, Conjecture \ref{conjv} plus the Bloch--Beilinson conjectures lead in particular to the following:

\begin{conjecture}\label{conj3} Let $X$ be a hyperk\"ahler variety of dimension $4$. Assume that $X$ admits a Lagrangian fibration with general fibre $A$. Then
\[ \ima \bigl( A^2(X)\xrightarrow{\cdot A} A^4(X)\bigr)=\QQ[c_4(T_X)]\ .\]
\end{conjecture}
(For a more general conjecture, which is more awkward to state, cf. \cite[Conjecture 1.3]{mmj}.)
  
 The goal of this note is to study the conjectural injectivity property (as outlined by Conjectures \ref{conjbv}, \ref{conjv} and \ref{conj3}) for some classical examples of Lagrangian fibrations.

%  Another conjecture concerns more generally $m$--dimensional subvarieties $Z\subset X$ that are {\em isotropic\/}, i.e. such that
%  \[  \cup Z\colon\ \ \ H^{2,0}(X)\ \to\ H^{m+2,m}(X) \]
%  is the zero map (isotropic subvarieties of dimension $m={1\over 2}\dim X$ are exactly the (possibly singular) Lagrangian subvarieties, cf. 
%  \cite[Proof of Lemma 1.4]{V17}).
%    Since $H^{\ast,0}(X)$ is generated by $H^{2,0}(X)$, we have that
%   \[ \cup Z\colon\ \ \ H^{j,0}(X)\ \to\ H^{m+j,m}(X) \]
%  is the zero map for all $j>0$.   
%  Since conjecturally, the piece $A^j_{[j]}(X)$ is determined by $H^{j,0}(X)$, and the piece $A^{m+j}_{[j]}(X)$ is determined by $H^{2m+j}(X)$, we arrive at the following conjecture:
%  
%  \begin{conjecture}\label{conjco} Let $X$ be a hyperk\"ahler variety of dimension $2m$. Let $Z\subset X$ be a codimension $m$ subvariety that is isotropic.
%   Then the maps
%    \[ \begin{split}    &A^j_{[j]}(X) \ \xrightarrow{\cdot Z}\  A^{m+j}(X)   \ \to\ A^{m+j}_{[j]}(X) \ ,\\
%                                  &A^m_{[j]}(X) \ \xrightarrow{\cdot Z}\  A^{2m}(X)   \ \to\ A^{2m}_{[j]}(X) \ ,\\
%       \end{split}    \]
%   are zero for all $j>0$. (Here, the right arrows are projection to the piece $A^{\ast}_{[j]}(X)$.)
%     \end{conjecture}
%  
  
  The first result is as follows:
  
   \begin{nonumbering}[=Theorem \ref{main1}] Let $X$ be a hyperk\"ahler fourfold, and assume that $X$ admits a Lagrangian fibration $\pi$ which is a compactified Jacobian of a family of curves. Let $A$ be a general fibre of $\pi$. Then
   \[ \ima \bigl( A^2(X)\xrightarrow{\cdot A} A^4(X)\bigr)=\QQ[c_4(T_X)]\ .\]
     \end{nonumbering}  
     
 %  Theorem \ref{main1} answers a particular case of a more general conjecture concerning isotropic subvarieties of hyperk\"ahler 
 %  varieties \cite[Conjecture 1.3]{mmj}. 
 To prove Theorem \ref{main1}, thanks to Markushevich \cite{Mark} one is reduced to fibrations arising from hyperplane sections of a genus $2$ K3 surface (such fibrations are cited as examples of Mukai flops in the introduction of Mukai's beautiful paper \cite[Example 0.6]{Mu}). Then, we exploit the existence of a {\em multiplicative Chow--K\"unneth decomposition\/} \cite{SV}, combined with results concerning the {\em Franchetta property\/} for families of Hilbert powers of low degree K3 surfaces \cite{FLV}.
     
 The second result is about the six-dimensional Lagrangian fibration $\pi\colon J_1\to \PP^3$, where $J_1$ is the compactified Jacobian of genus $3$ curves arising as hyperplane sections of a general quartic K3 surface $S$. This is another example in the introduction of Mukai's foundational paper \cite[Example 0.8]{Mu}, where it is shown that the flop of $J_1$ along a certain codimension $2$ subvariety $P\subset J_1$ is isomorphic to a moduli space of sheaves on $S$.
  
 \begin{nonumbering}[=Theorem \ref{main2}] Let $h_1\in A^1(J_1)$ be the polarization class, let $h_2\in A^1(J_1)$ be $\pi^\ast(d)$ where $d\subset\PP^3$ is a hyperplane class, and let
 $P\subset J_1$ be as above.
%let $A\in A^3(J_1)$ be a general fibre of the fibration $J_1\to\PP^3$. 
Let $R^\ast(J_1)$ be the $\QQ$-subalgebra 
   \[ R^\ast(J_1):=  \langle  h_1, h_2, P, c_j(T_{J_1})\rangle\ \ \ \ \subset\ A^\ast(J_1)\ .\]
   The cycle class map induces an injection
   \[ R^\ast(J_1)\ \hookrightarrow\ H^\ast(J_1,\QQ)\ .\]
 \end{nonumbering}

  In particular, Conjecture \ref{conjbv} is true for the very general sixfold $J_1$.
  Theorem \ref{main2} is also in agreement with Conjecture \ref{conjv}, because $P$ (being a $\PP^2$-bundle over $S$) has codimension $2$ and is swept out by constant cycle surfaces. 
 In proving Theorem \ref{main2}, we rely on (a sharpening of) a recent result of B\"ulles \cite{Bul}, combined with results on the Franchetta property from \cite{FLV}.
 
 For Lagrangian fibrations of higher dimension (such as the tenfolds of \cite{LSV} or \cite{V18}), the argument of the present note quickly runs into problems: this is because the Franchetta property is not  known outside of a few selected cases (e.g., for the tenfolds of \cite{V18}, in view of \cite{LPZ} one would need the Franchetta property for the 5th relative power of cubic fourfolds; this is currently unknown and perhaps not even true).
    
\vskip0.6cm

\begin{convention} In this note, the word {\sl variety\/} will refer to a reduced irreducible scheme of finite type over $\C$. For a smooth variety $X$, we will denote by $A^j(X)$ the Chow group of codimension $j$ cycles on $X$ 
with $\QQ$-coefficients.
%for $X$ smooth of dimension $n$ the notations $A_jX$ and $A^{n-j}X$ will be used interchangeably. 
%In the rare cases we will have something to say about Chow groups with integral coefficients, we will indicate this by writing $A^\ast X_{\ZZ}$.
The notation 
%$A^j_{num}(X)$, 
$A^j_{hom}(X)$ will be used to indicate the subgroups of 
%numerically resp. 
homologically trivial cycles.

For a morphism between smooth varieties $f\colon X\to Y$, we will write $\Gamma_f\in A^\ast(X\times Y)$ for the graph of $f$, and ${}^t \Gamma_f\in A^\ast(Y\times X)$ for the transpose correspondence.

The contravariant category of Chow motives (i.e., pure motives with respect to rational equivalence as in \cite{Sc}, \cite{MNP}) will be denoted $\MM_{\rm rat}$. 
%We will write $H^\ast(X)=H^\ast(X,\QQ)$ for singular cohomology with $\QQ$--coefficients.
\end{convention}

\section{Preliminaries}

\subsection{B\"ulles' result revisited}

The following theorem is a slight sharpening of a result of B\"ulles \cite{Bul}:

\begin{theorem}
%[B\"ulles \cite{Bul}]
\label{bul} Let $S$ be a projective K3 surface or an abelian surface, and $\alpha\in Br(S)$ a Brauer class. Let $M$ be a smooth projective moduli space
of Gieseker stable $\alpha$-twisted sheaves on $S$, of dimension $\dim M=2m$. There is an inclusion as direct summand
  \[   h(M)\ \hookrightarrow\ \bigoplus_{i=1}^r h(S^{k_i})(\ell_i)\ \ \ \hbox{in}\ \MM_{\rm rat}\ ,\]
  where $\ell_i\in\Z$ and $1\le k_i\le m$.
  \end{theorem}
  
  \begin{proof} We follow B\"ulles' proof, with a slight twist to get a better bound on the integers $k_i$ (in \cite[Theorem 0.1]{Bul}, the $k_i$ are bounded by $2m$).
  Let
    \[  [\hbox{Ext}^!_\pi]:=\sum_i (-1)^i [\hbox{Ext}^i_\pi(\EE,\FF)]   \ \ \in\ K_0(M\times M) \]
    be as in \cite[Proof of Theorem 0.1]{Bul}. Then (as explained in loc. cit.) a result of Markman's gives the equality
    \begin{equation}\label{mm} \Delta_M = c_{2m}(- [\hbox{Ext}^!_\pi]) \ \ \ \hbox{in}\   A^{2m}(M\times M)\ .\end{equation}
    
As in loc. cit., we consider the two-sided ideal in the ring of correspondences
  \[     I:=  \bigcup_{k\ge 1} I_k \ \ \  \subset \ A^{\ast}(M\times M)\ ,\] 
  where $I_k$ is defined as
  \[ I_k:= \langle \beta\circ\alpha\ \big\vert\ \alpha\in A^\ast(M\times S^\ell)\ ,\ \beta\in A^\ast( S^\ell\times M)\ ,\  1\le \ell\le k \rangle\ .\]
  B\"ulles shows \cite[Proof of Theorem 0.1]{Bul} that $I$ is closed under intersection product, and more precisely that $I_k\cdot I_\ell\subset I_{k+\ell}$.
  In addition, let us state a lemma:
  
  \begin{lemma}\label{div} Let $\gamma\in I_k$ (for some $k\ge 1$), and $\delta\in A^1(M\times M)$. Then
    \[ \gamma\cdot\delta\ \ \in\ I_k\ .\]
    \end{lemma}
    
    \begin{proof} Since the irregularity $q(M)=0$, every divisor $\delta\in A^1(M\times M)$ can be written as a sum of pullbacks $(p_i)^\ast(D_i)$ under the two projections
    $p_i\colon M\times M\to M$. We may thus suppose $\delta$ is of the form $D\times M$ or $M\times D$, where $D\subset M$ is an irreducible reduced divisor. Let $\iota\colon D\to M$ denote the inclusion morphism. We have
     \[ \begin{split}  \gamma\cdot\delta=\gamma\cdot (D\times M) &=  (\iota\times\ide)_\ast(\iota\times\ide)^\ast(\gamma)=    (\Gamma_\iota\times\Delta_M)_\ast ({}^t \Gamma_\iota\times\Delta_M)_\ast  (\gamma)\\
                         &=  \Gamma_\iota\circ   {}^t \Gamma_\iota\circ \gamma  \ \ \ \ \hbox{in}\ A^\ast(M\times M)\ ,\\
                         \end{split}\]
                         where the last equality is by virtue of Lieberman's lemma \cite[Lemma 3.3]{V3}, \cite[Proposition 2.1.3]{MNP}. Similarly, in case $\delta$ is of the form $M\times D$, we find that
                         \[ \gamma\cdot\delta=\gamma\cdot (M\times D)= \gamma\circ  \Gamma_\iota\circ   {}^t \Gamma_\iota\ \ \ \hbox{in}\ A^\ast(M\times M)\ .\]
             In both cases, it follows that $\gamma\in I_k$ implies that also $\gamma\cdot \delta\in I_k$.        
                          \end{proof}
  
  Let us write $c_n:=c_n(- [\hbox{Ext}^!_\pi])\in A^n(M\times M)$.  As shown by B\"ulles, we have
    \begin{equation}\label{rel} c_n=(-1)^{n-1}(n-1)!\, ch_n + p(c_1,\ldots,c_{n-1})\ \ \ \hbox{in}\ A^{n}(M\times M)\ \ \ \forall n\ge 1\ .\end{equation}
   Here $ch_n$ denotes the degree $n$ part of the Chern character $ch(- [\hbox{Ext}^!_\pi])\in A^\ast(M\times M)$, and $p$ is some weighted homogeneous polynomial of degree $n$. We have $ch_n\in I_1$ for all $n\ge 1$.
    In particular, for $n=2$ we find that
     \[ c_2 = {1\over 2} c_1^2 - ch_2\ \ \ \hbox{in}\ A^{2}(M\times M)\ .\]
    The class $c_1=ch_1$ is in $I_1$ and so (using lemma \ref{div}) $c_1^2$ is also in $I_1$. It follows that 
    \[ c_2\ \ \in\ I_1\ .\]
    Likewise, $c_3$ can be expressed in terms of $ch_3\in I_1$ and $c_1^3\in I_1$ and $c_1\cdot c_2\in I_1$, and so (again using lemma \ref{div}) we see that
     \[ c_3\ \ \in I_1\ .\]
   
   We now make the claim that
     %  \begin{equation}\label{this} c_n\ \ \in\ I_{\lfloor {n\over 2}\rfloor}\ \ \ \forall n\ge 2\ .\end{equation}
     \begin{equation}\label{thispol}   p(c_1,\ldots,c_{n-1})   \ \ \in\ I_{\lfloor {n\over 2}\rfloor}\ \end{equation}
     for any weighted homogeneous polynomial $p$ of degree $n\ge 2$. Let us prove this claim by induction.
     From what we have just checked, it is clear that the claim is true for $n=2, 3$. Let us now suppose $n\ge 4$. The polynomial $p$ can be decomposed
    \[ p(c_1,\ldots,c_{n-1})= \lambda c_1\cdot c_{n-1}+\mu c_1^2c_{n-2} +\nu c_2 c_{n-2} + c_1^2 q(c_1,\ldots,c_{n-3})+ c_2 r(c_1,\ldots,c_{n-3})\ ,\]
    where $\lambda, \mu,\nu\in\QQ$ and $q$ and $r$ are weighted homogeneous polynomials of degree $n-2$. By the induction hypothesis combined with (\ref{rel}), we know that $c_{n-1}$ and $c_{n-2}$ are in $I_{\lfloor {n\over 2}\rfloor}$; using lemma \ref{div} this implies that the terms $ c_1\cdot c_{n-1}$ and $ c_1^2c_{n-2}$ are in $I_{\lfloor {n\over 2}\rfloor}$. By the induction hypothesis, $c_{n-2}\in I_{\lfloor {n-2\over 2}\rfloor}$; since $c_2\in I_1$ this gives $ c_2 c_{n-2}\in I_{\lfloor {n\over 2}\rfloor}$. Again by the induction hypothesis, the polynomials $q$ and $r$ are in $ I_{\lfloor {n-2\over 2}\rfloor}$.
    Using lemma \ref{div}, it follows that $ c_1^2 q(c_1,\ldots,c_{n-3})\in  I_{\lfloor {n-2\over 2}\rfloor}$. Using the fact that $c_2\in I_1$, it follows that $c_2 r(c_1,\ldots,c_{n-3})\in   I_{\lfloor {n\over 2}\rfloor}$. Altogether, this proves the claim (\ref{thispol}).     
    
 Claim (\ref{thispol}), combined with relation (\ref{rel}) and the fact that $ch_n\in I_1$, implies that
    \[ c_n\ \ \in\ I_{\lfloor {n\over 2}\rfloor}\ \ \ \forall n\ge 2\ .\]
   In view of equality (\ref{mm}), it follows that
   \[ \Delta_M\ \ \in\ I_m\cap A^{2m}(M\times M)\ ,\]
   which proves the theorem.   
 \end{proof}
 
 \begin{remark} As noted by B\"ulles \cite[Remark 2.1]{Bul}, Theorem \ref{bul} is also valid for moduli spaces of $\sigma$-stable objects on a K3 surface or abelian surface, where $\sigma$ is a generic stability condition. For instance, Theorem \ref{bul} applies to Ouchi's eightfolds \cite{Ouch} and to the Laza--Sacc\`{a}--Voisin tenfolds \cite{LSV}.
 
 In \cite{FLV2}, we use Theorem \ref{bul} to prove the generalized Franchetta conjecture for Lehn--Lehn--Sorger--van Straten eightfolds.
 \end{remark}

  \subsection{The Franchetta property}
  
  \begin{definition} Let $\pi\colon\XX\to B$ be a smooth projective family of varieties, and let us write $X_b:=\pi^{-1}(b)$ for a fibre. We say that the family $\pi\colon\XX\to B$ has the {\em Franchetta property\/} if for any $\Gamma\in A^\ast(\XX)$ there is equivalence
  \[ \begin{split}\Gamma\vert_{X_b}=0\ \ \hbox{in}\ H^\ast(X_b)\ \ \ \hbox{for\ $b\in B$\ very\ general}\ \ \iff\ \  \Gamma\vert_{X_b}=0\ \ \hbox{in}\ A^\ast(X_b)  \ \ \ \hbox{for}&\\
  \hbox{$b\in B$\ very\ general}&\ .\\
  \end{split}\]
  \end{definition}
  
  \begin{remark} In view of \cite[Lemma 3.2]{Vo}, the vanishing $ \Gamma\vert_{X_b}=0$ in $A^\ast(X_b)$ for $b\in B$ very general is equivalent to the vanishing $ \Gamma\vert_{X_b}=0$ in $A^\ast(X_b)$ for all $b\in B$.  
  \end{remark}

  \begin{notation}\label{not} Let $\PP$ denote weighted projective space $\PP(1^3,3)$. Let
  $\Ss_{g2}\to B_{g2}$ denote the universal family of K3 surfaces of genus $2$, where 
    \[B_{g2}\subset \PP H^0(\PP,\OO_\PP(6))\] 
    is the Zariski open parametrizing smooth sections.
  
  Let $\Ss_{g3}\to B_{g3}$ denote the universal family of K3 surfaces of genus $3$, where 
    \[B_{g3}\subset \PP H^0(\PP^3,\OO_{\PP^3}(4))\] 
    is the Zariski open parametrizing smooth sections.  
    \end{notation}
  
  \begin{notation} For any family $\XX\to B$ and $m\in\NN$, we write $\XX^{m/B}:=\XX\times_B\cdots\times_B \XX$ for the $m$-fold fibre product.
  \end{notation}

  \begin{theorem}[\cite{FLV}]\label{gfc} The families $\Ss_{g2}^{m/B_{g2}}\to B_{g2}$, $m\le 3$ and $\Ss_{g3}^{m/B_{g3}}\to B_{g3}$, $m\le 5$
  have the Franchetta property.   
  \end{theorem}
    
   \begin{proof} This is (part of) \cite[Theorem 1.5]{FLV}.   
       \end{proof} 
    
   \section{Multiplicative Chow--K\"unneth decomposition}

	\begin{definition}[Murre \cite{Mur}]\label{ck} Let $X$ be a smooth projective
		variety of dimension $n$. We say that $X$ has a 
		{\em CK  decomposition\/} if there exists a decomposition of the
		diagonal
		\[ \Delta_X= \pi^0_X+ \pi^1_X+\cdots +\pi^{2n}_X\ \ \ \hbox{in}\
		A^n(X\times X)\ ,\]
		such that the $\pi^i_X$ are mutually orthogonal idempotents and
		$(\pi^i_X)_\ast H^\ast(X)= H^i(X)$.
		Given a CK decomposition for $X$, we set 
		$$A^i_{(j)}(X) := (\pi_X^{2i-j})_\ast A^i(X).$$
		The CK decomposition is said to be {\em self-dual\/} if
		\[ \pi^i_X = {}^t \pi^{2n-i}_X\ \ \ \hbox{in}\ A^n(X\times X)\ \ \ \forall
		i\ .\]
		(Here ${}^t \pi$ denotes the transpose of a cycle $\pi$.)
		
		  (NB: ``CK decomposition'' is short-hand for ``Chow--K\"unneth
		decomposition''.)
	\end{definition}
	
	\begin{remark} \label{R:Murre} The existence of a Chow--K\"unneth decomposition
		for any smooth projective variety is part of Murre's conjectures \cite{Mur},
		\cite{MNP}. 
		It is expected that for any $X$ with a CK
		decomposition, one has
		\begin{equation*}\label{hope} A^i_{(j)}(X)\stackrel{??}{=}0\ \ \ \hbox{for}\
		j<0\ ,\ \ \ A^i_{(0)}(X)\cap A^i_{num}(X)\stackrel{??}{=}0.
		\end{equation*}
		These are Murre's conjectures B and D, respectively.
	\end{remark}

	\begin{definition}[Definition 8.1 in \cite{SV}]\label{mck} Let $X$ be a
		smooth
		projective variety of dimension $n$. Let $\Delta_X^{sm}\in A^{2n}(X\times
		X\times X)$ be the class of the small diagonal
		\[ \Delta_X^{sm}:=\bigl\{ (x,x,x) : x\in X\bigr\}\ \subset\ X\times
		X\times X\ .\]
		A CK decomposition $\{\pi^i_X\}$ of $X$ is {\em multiplicative\/}
		if it satisfies
		\[ \pi^k_X\circ \Delta_X^{sm}\circ (\pi^i_X\otimes \pi^j_X)=0\ \ \ \hbox{in}\
		A^{2n}(X\times X\times X)\ \ \ \hbox{for\ all\ }i+j\not=k\ .\]
		In that case,
		\[ A^i_{(j)}(X):= (\pi_X^{2i-j})_\ast A^i(X)\]
		defines a bigraded ring structure on the Chow ring\,; that is, the
		intersection product has the property that 
		\[  \ima \Bigl(A^i_{(j)}(X)\otimes A^{i^\prime}_{(j^\prime)}(X)
		\xrightarrow{\cdot} A^{i+i^\prime}(X)\Bigr)\ \subseteq\ 
		A^{i+i^\prime}_{(j+j^\prime)}(X)\ .\]
		
		(For brevity, we will write {\em MCK decomposition\/} for ``multiplicative Chow--K\"unneth decomposition''.)
	\end{definition}
	
	\begin{remark}
	The property of having an MCK decomposition is
	severely restrictive, and is closely related to Beauville's ``splitting
	property'' conjecture \cite{Beau3}. Examples of varieties admitting an MCK decomposition include hyperelliptic curves, K3 surfaces, abelian varieties, cubic hypersurfaces.
	For more ample discussion and more examples, we refer to \cite[Chapter 8]{SV}, as well as \cite{V6}, \cite{SV2}, \cite{FTV}, \cite{FV}, \cite{LV}, \cite{FLV}.
	\end{remark}

There are the following useful general results:	

\begin{theorem}[Shen--Vial \cite{SV}]\label{square} Let $X$ be a hyperk\"ahler fourfold that is birational to a Hilbert square $S^{[2]}$ where $S$ is a K3 surface. Then $X$ has an MCK decomposition.
\end{theorem}

\begin{proof} The statement for $S^{[2]}$ is \cite[Theorem 13.4]{SV} (a more general result is \cite[Theorem 1]{V6}). The statement for $X$ then follows by applying the result of Rie\ss\  \cite{Rie} (as duly noted in \cite[Introduction]{V6}).
\end{proof}

\begin{proposition}[Shen--Vial \cite{SV}]\label{product} Let $M,N$ be smooth projective varieties that have an MCK decomposition. Then the product $M\times N$ has an MCK decomposition.
\end{proposition}

\begin{proof} This is \cite[Theorem 8.6]{SV}, which shows more precisely that the {\em product CK decomposition\/}
  \[ \pi^i_{M\times N}:= \sum_{k+\ell=i} \pi^k_M\times \pi^\ell_N\ \ \ \in A^{\dim M+\dim N}\bigl((M\times N)\times (M\times N)\bigr) \]
  is multiplicative.
\end{proof}

\begin{theorem}[Shen--Vial \cite{SV2}]\label{blowup} Let $M$ be a smooth projective variety, and let $f\colon\wt{M}\to M$ be the blow--up with center a smooth closed subvariety
$N\subset M$. Assume that
\begin{enumerate}

\item $M$ and $N$ have a self-dual MCK decomposition;

\item the Chern classes of the normal bundle $\NNN_{N/M}$ are in $A^\ast_{(0)}(N)$;

\item the graph of the inclusion morphism $N\to M$ is in $A^\ast_{(0)}(N\times M)$.

%\item the Chern classes $c_j(T_M)$ are in $A^\ast_{(0)}(M)$.

\end{enumerate}
Then $\wt{M}$ has a self-dual MCK decomposition, and 
% the graph $\Gamma_f$ is in $A^\ast_{(0)}(  \wt{M}\times M)$.
\[  f^\ast A^\ast_{(j)}(M)\ \subset\ A^\ast_{(j)}(\wt{M})\ ,\ \ \      f_\ast A^\ast_{(j)}(\wt{M})\ \subset\ A^\ast_{(j)}({M})\ .\]
\end{theorem}

\begin{proof} This is \cite[Proposition 2.4]{SV2}.
\end{proof}

\section{Examples in dimension 4}

 \begin{theorem}\label{main1} Let $X$ be a hyperk\"ahler fourfold, and assume that $X$ admits a Lagrangian fibration $\pi\colon X\to B$ which is a compactified Jacobian of a family of curves. Let $A$ be a general fibre of $\pi$. Then
   \[ \ima \bigl( A^2(X)\xrightarrow{\cdot A} A^4(X)\bigr)=\QQ[c_4(T_X)]\ .\]
     \end{theorem}
 
 \begin{proof} Thanks to a result of Markushevich \cite[Theorem 1.1]{Mark}, we know that $B\cong\PP^2$ and $X\cong J_0$, where $J_0$ is the compactified Jacobian of the genus $2$ curves arising as hyperplane sections of a genus $2$ K3 surface $S$. The fibration $\pi\colon J_0\to\PP^2$ occurs in \cite[Example 0.6]{Mu}, where it is shown that 
 %$J_0$ is related to the Hilbert square $S^{[2]}$ via 
 there is a birational map
     \[  J_0\ \dashrightarrow\ S^{[2]} \]
  which is a Mukai flop. Precisely, $S^{[2]}$ contains a subvariety $P\cong\PP^2$ (defined as the pairs of points in $S$ that are in the same fibre of the double cover $S\to\PP^2$). There are birational transformations
  \[  J_0\ \xleftarrow{r}\ \wt{J_0}\ \xrightarrow{s}\ S^{[2]}\ ,\]
  where $s$ is the blow-up with center $P\subset S^{[2]}$, and $r$ is the blow-down of the exceptional divisor of $s$ onto a closed subvariety
  $P^\prime\subset J_0$.

%  We know from Rie\ss's work \cite{Rie} that there exists a correspondence $\Gamma$ inducing a ring isomorphism
%   \[ \Gamma_\ast\colon\ \ A^\ast(J_0)\ \cong\ A^\ast(S^{[2]}\ .\]
%   (Since $f$ is a Mukai flop, this does not need \cite[Theorem 3.2]{Rie}; this follows from the more explicit argument \cite[Proposition 6.6]{Rie}.)
%  Since $\Gamma$ also induces an isomorphism of Chow motives, it follows 
%We know from theorems \ref{square} and \ref{blowup} that $\wt{J_0}$ has an MCK decomposition, induced by the MCK decomposition of $S^{[2]}$ constructed in \cite{SV}. 
%Let us write
%  $A^\ast_{(\ast)}(\wt{J_0})$ for the bigrading of the Chow ring induced by this MCK decomposition.
  
  The theorem will follow by combining the following 4 claims:
  
  \begin{claim}\label{c0} The variety $\wt{J_0}$ has an MCK decomposition, and this induces a splitting $A^2(\wt{J_0})=A^2_{(0)}(\wt{J_0})\oplus A^2_{(2)}(\wt{J_0})$.
  \end{claim}
  
  \begin{claim}\label{c1} Let $A$ be a fibre of the Lagrangian fibration $\pi\colon J_0\to\PP^2$. Then
    \[ r^\ast(A)\  \ \in\ A^2_{(0)}(\wt{J_0})\ .\]
    \end{claim}
    
    \begin{claim}\label{c2} Let $A$ be a general fibre of $\pi\colon J_0\to\PP^2$. The map
      \[ A^2_{(2)}(\wt{J_0})\ \xrightarrow{\cdot r^\ast(A)}\ A^4(\wt{J_0}) \]
      is zero.
      \end{claim}    
      
      \begin{claim}\label{c3} One has 
        \[r^\ast c_4(T_{J_0})\in A^4_{(0)}(\wt{J_0}) \ .\]
        \end{claim}
        
    Let us show that these claims imply the theorem; Since $A^2(\wt{J_0})=A^2_{(2)}(\wt{J_0})\oplus A^2_{(0)}(\wt{J_0})$, we have
     \[ \ima \bigl( A^2(\wt{J_0})\xrightarrow{\cdot r^\ast(A)} A^4(\wt{J_0})\bigr) =  A^2_{(2)}(\wt{J_0})\cdot r^\ast(A) + A^2_{(0)}(\wt{J_0})\cdot r^\ast(A)\ .\]     
        Using Claim \ref{c2}, this reduces to
    \[ \ima \bigl( A^2(\wt{J_0})\xrightarrow{\cdot r^\ast(A)} A^4(\wt{J_0})\bigr) =    A^2_{(0)}(\wt{J_0})\cdot r^\ast(A)\ .\]     
    Using Claim \ref{c1}, we see that $A^2_{(0)}(\wt{J_0})\cdot r^\ast(A)$ is contained in $A^4_{(0)}(\wt{J_0})\cong\QQ$. 
    %Because the intersection of $r^\ast(A)$ with 2 ample divisors has strictly positive degree, there is equality 
   % \[  A^2(\wt{J_0})\cdot r^\ast(A)=  A^2_{(0)}(\wt{J_0})\cdot r^\ast(A)\cong A^4_{(0)}(\wt{J_0})\cong\QQ  \ .\]
    Claim \ref{c3},  plus the fact that $c_4(T_{J_0})$ has strictly positive degree, then implies that
    \[  A^2(\wt{J_0})\cdot r^\ast(A) \ \ \in\    \QQ[ r^\ast c_4(T_{J_0})]\ .\]
    Pushing forward to $J_0$, this gives an inclusion
    \[ A^2({J_0})\cdot A \ \ \in\    \QQ[  c_4(T_{J_0})]\ \ \subset\ A^4(J_0)\ .\]     
   Since the left-hand side is one-dimensional (the intersection of $A$ with 2 ample divisors has strictly positive degree), this inclusion is an equality, proving the theorem.    
    
   It remains to prove the claims. 
   To prove Claim \ref{c0}, we use Theorem \ref{blowup} with $M=S^{[2]}$ and $N=P\cong\PP^2$. Points (1) and (2) are clearly satisfied. For point (3), we note that $S^{[2]}$ and $P$ have a ``universal MCK decomposition'', i.e. there exist
    \[ \pi^i_{\Ss^{[2]}\times_B \PPP}\ \ \ \in\ A^6(\Ss^{[2]}\times_B \PPP)\ \ , \ \ i=0,\ldots,12\ ,\]
    such that for each $b\in B$ the restriction
    \[ \pi^i_{(S_b)^{[2]}\times P_b }:=  \pi^i_{\Ss^{[2]}\times_B \PPP}\vert_b\ \ \ \in\ A^6((S_b)^{[2]}\times P_b   )\ \]
    defines an MCK decomposition for  $(S_b)^{[2]}\times P_b$. Let $\iota\colon \PPP\to \Ss^{[2]}$ denote the inclusion morphism, and $\iota_b\colon P_b\to (S_b)^{[2]}$ the restriction to a fibre. For any $k\not=8$, we have that
    \[  (\pi^k_{(S_b)^{[2]}\times P_b })_\ast (\Gamma_{\iota_b}) = \Bigl(  (\pi^k_{\Ss^{[2]}\times_B \PPP})_\ast (\Gamma_\iota)\Bigr)\vert_b\ \ \ \in\ A^4((S_b)^{[2]}\times P_b) \]
    is homologically trivial. Theorem \ref{gfc} then implies that it is rationally trivial, and so
    \[ \Gamma_{\iota_b}=    (\pi^8_{(S_b)^{[2]}\times P_b })_\ast (\Gamma_{\iota_b}) \ \ \ \hbox{in}\ A^4((S_b)^{[2]}\times P_b)\ \ \ \forall b\in B \ .\]
    We have now checked that the conditions of Theorem \ref{blowup} are satisfied, and so $\wt{J_0}$ has an MCK decomposition. The ``blow-up'' isomorphism $A^2(\wt{J_0})\cong A^2(S^{[2]})\oplus A^1(P)$ is homogeneous with respect to the lower grading. Since $A^2_{(j)}(S^{[2]})=0$ for $j\not\in\{0,2\}$ and
    $A^1(P)=A^1_{(0)}(P)$, this shows the second part of claim \ref{c0}.
   
   Claim \ref{c1} is elementary: writing $A=\pi^{-1}(x)$ where $x\in\PP^2$, we see that
   $r^\ast(A)=r^\ast\pi^\ast(x)$ in $A^2(\wt{J_0})$ is an intersection of divisors, which proves the claim.
   To prove the other 2 claims, we consider things familywise. That is, we let $\Ss\to B$ denote the universal family of genus $2$ K3 surfaces as in notation \ref{not}, and we write $\Ss^{[2]}\to B$ for the universal family of Hilbert squares of genus $2$ K3 surfaces. There are morphisms of $B$-schemes
    \[  \JJ\ \xleftarrow{r}\ \wt{\JJ}\ \xrightarrow{s}\ \Ss^{[2]}\ ,\]   
    such that restriction to a fibre gives the Mukai flop mentioned above. The morphism $r$ is the blow-up with center $\PPP^\prime$ (which is a $\PP^2$-bundle over $B$), and the morphism $s$ is the blow-up with center $\PPP$ (which is again a $\PP^2$-bundle over $B$).
   We now establish the following result:
   
   \begin{proposition}\label{gfc4} Let $\wt{\JJ}\to B$ be as above. The families $\wt{\JJ}^{}\to B$ 
   and $\wt{\JJ}\times_B \Ss\to B$ have the Franchetta property.
   \end{proposition}
   
   \begin{proof}   
   %  To prove the proposition, it will suffice to prove the Franchetta property for the families $\wt{\JJ}^{}\to B$ and $\wt{\JJ}\times_B \Ss\to B$ .
    For the first family, one notes that there is a commutative diagram
      \[
   \begin{array}[c]{ccc}
      A^i(\wt{\JJ}) & \xrightarrow{}& A^i(\Ss^{[2]})\oplus A^{i-1}(\PPP)\\
       &&\\
       \downarrow&&\downarrow\\
       &&\\
      A^i(\wt{J}_b) & \xrightarrow{\cong}& A^i((S_b)^{[2]})\oplus A^{i-1}(P_b)\ .\\   
      \end{array}\]     
      (Here, we write $\wt{J}_b, S_b, P_b$ for the fibre over $b\in B$ of the family $\wt{\JJ}$ resp. $\Ss$ resp. $\PPP$.)
      The family $\Ss^{[2]}\to B$ has the Franchetta property (theorem \ref{gfc}), and $P_b\cong\PP^2$ so the family $\PPP\to B$ trivially has the Franchetta property. This settles the Franchetta property for $\wt{\JJ}\to B$.
      
      For the second family, there is a similar commutative diagram
      \[
   \begin{array}[c]{ccc}
      A^i(\wt{\JJ}\times_B \Ss) & \xrightarrow{}& A^i(\Ss^{[2]}\times_B \Ss)\oplus A^{i-1}(\PPP\times_B \Ss)\\
       &&\\
       \downarrow&&\downarrow\\
       &&\\
      A^i(\wt{J}_b\times S_b) & \xrightarrow{\cong}& A^i((S_b)^{[2]}\times S_b)\oplus A^{i-1}(P_b\times S_b)\ .\\   
      \end{array}\]     
    The family $\Ss^{[2]}\times\Ss\to B$ has the Franchetta property (theorem \ref{gfc}, or more exactly \cite[Theorem 1.5]{FLV}), and so does the family $\PPP\times_B \Ss\to B$ (using the projective bundle formula, one reduces to $\Ss\to B$). This settles the Franchetta property for $\wt{\JJ}\times_B \Ss\to B$.
       \end{proof}
   
   Let us now prove the two remaining claims. We will rely on the existence of an MCK decomposition that is {\em generically defined\/} for the family $\wt{\JJ}\to B$, in the following sense:
         
      \begin{lemma}\label{umck} Let $\wt{\JJ}\to B$ be as above. There exist
    \[ \pi^i_{\wt{\JJ}}\ \ \ \in\ A^4(\wt{\JJ}\times_B \wt{\JJ})\ \ , \ \ i=0,\ldots,8\ ,\]
    such that for each $b\in B$ the restriction
    \[ \pi^i_{\wt{J_b}}:=  \pi^i_{\wt{\JJ}}\vert_b\ \ \ \in\ A^4(\wt{J}_b\times \wt{J}_b)\ \]
    defines an MCK decomposition for $\wt{J}_b$.
    \end{lemma}
    
    \begin{proof} The Hilbert squares $(S_b)^{[2]}$ have an MCK decomposition that exists universally (this is just because the ``distinguished $0$-cycle'' of \cite{BV} exists universally).
  Looking at the argument of Theorem \ref{blowup} (i.e. the proof of \cite[Proposition 2.4]{SV2}, one sees that the induced MCK decomposition for the blow-up $\wt{J_b}$ exists universally as well.
    \end{proof}
    
To prove Claim \ref{c3}, we observe that 
   \[ (r_b)^\ast  c_4(T_{J_b})=\bigl( r^\ast c_4(T_{\JJ/B})\bigr)\vert_b\ \ \ \in A^4(\wt{J}_b)\  \]   
   is universally defined. This forces $(r_b)^\ast  c_4(T_{J_b})$ to lie in $ A^4_{(0)}(\wt{J}_b) $: for any $k\not=8$, we have that 
   \[  (\pi^k_{\wt{J}_b})_\ast  \bigl((r_b)^\ast  c_4(T_{J_b})\bigr) =  \Bigl( (\pi^k_{\wt{\JJ}})_\ast (r^\ast  c_4(T_{\JJ/B})) \Bigr)\vert_b  \ \ \ \in A^4(\wt{J}_b) \]
   is homologically trivial, for all $b\in B$. In view of Proposition \ref{gfc4}, this implies 
    \[  (\pi^k_{\wt{J}_b})_\ast  \bigl((r_b)^\ast  c_4(T_{J_b})\bigr)=0\ \ \ \hbox{in}\ A^4(\wt{J}_b)\ \ \ \forall b\in B\ ,\ \ \ \forall k\not=8\ ,\]
    proving claim \ref{c3}.   
    
    To prove Claim \ref{c2}, let $\AAA\subset \JJ$ denote a general fibre of $\pi\colon\JJ\to \cup_{b\in B} \vert \OO_{J_b}(1)\vert$, and let
    $\wt{\AAA}$ denote a general fibre of $\pi\circ r$. Let us write $\tau\colon \wt{\AAA}\to\wt{\JJ}$ for the inclusion. We are interested in the correspondence
    \[ \Gamma_b:=      \Gamma_{\tau_b} \circ {}^t \Gamma_{\tau_b}   \circ \pi^2_{\wt{J}_b} \ \ \ \in\ A^6(\wt{J}_b\times \wt{J}_b)\ ,\]
    which by construction is such that
    \[ (\Gamma_b)_\ast A^2(\wt{J}_b) = A^2_{(2)}(\wt{J}_b)\cdot (r_b)^\ast (A_b)\ .\]
    The correspondence $\Gamma_b$ is universally defined, i.e. there exists $\Gamma\in A^6(\wt{\JJ}\times_B \wt{\JJ})$ such that 
      \[ \Gamma_b = \Gamma\vert_b\ \ \ \in\ A^6(\wt{J}_b\times \wt{J}_b)\ \ \ \forall b\in B\ .\]
    Since $A_b\subset J_b$ is Lagrangian, the cup product of $A_b$ with $H^{2,0}(J_b)$ is zero. By a standard Hodge theory argument, this means that the cup product of $A_b$ with the
    transcendental cohomology $H^2_{tr}(J_b)$ is also zero. Since $H^2_{tr}$ is a birational invariant, the same holds on $\wt{J_b}$, and so the map
    \[ H^2(\wt{J_b})\ \xrightarrow{\cdot\wt{A_b}}\ H^6(\wt{J_b}) \]
    is the same as the map
    \[ N^1 H^2(\wt{J_b})\ \xrightarrow{\cdot\wt{A_b}}\ N^3 H^6(\wt{J_b}) \]
    (where $N^i H^{2i}()$ denotes the algebraic classes in cohomology).
    It follows that there exist (for each $b\in B$) a finite union of curves $C_b\subset\wt{J}_b$ and a cycle $\gamma_b$ supported on $C_b\times C_b$ such that
     \[ \Gamma_b=\gamma_b\ \ \ \hbox{in}\ H^{12}(\wt{J}_b\times\wt{J}_b)\ .\]
     (Indeed, for $C_b$ one can take a basis of $N^3 H^6(\wt{J_b})$, and add curves forming a dual basis to $N^1 H^2(\wt{J_b})$.)
          Using Voisin's Hilbert schemes argument as in \cite[Proposition 3.7]{V0}, these fibrewise data can be spread out, i.e. there exist a finite union of codimension $3$ closed subvarieties $\CCC\subset\wt{\JJ}$ and a cycle $\gamma$ supported on $\CCC\times_B \CCC\subset \wt{\JJ}\times_B \wt{\JJ}$ with the property that
     \begin{equation}\label{homtriv}  (\Gamma -\gamma)\vert_b =0\ \ \ \hbox{in}\ H^{12}(\wt{J}_b\times\wt{J}_b)\ \ \ \forall b\in B\ .\end{equation}   
     At this point, we need another lemma:
     
     \begin{lemma}\label{ll} Set-up as above. There exist relative correspondences
    \[  \Theta_1\ ,\ \Theta_2\in A^{4}(\Ss\times_B \wt{\JJ})\ ,\ \ \ \Xi_1\ , \ \Xi_2\in  A^{2}( (\wt{\JJ}\times_B \Ss)  \]
    such that for each $b\in B$, the composition
    \[   A^{2}_{(2)}(\wt{J_b})\ \xrightarrow{\bigl((\Xi_1\vert_b)_\ast, (\Xi_2\vert_b)_\ast\bigr)}\
             A^2(S_b)\oplus  A^2(S_b)
              \xrightarrow{\bigl((\Theta_1+\Theta_2)\vert_b\bigr)_\ast}\ A^{2}(\wt{J_b}) \]
      is the identity.
      \end{lemma}
      
      \begin{proof} By virtue of Theorem \ref{blowup}, the isomorphism
      \[ A^2(\wt{J_b})\ \cong\ A^2\bigl((S_b)^{[2]}\bigr)\oplus A^1(P_b) \]
      respects the bigrading. Since $A^1_{(2)}(P_b)=0$, it follows that  
         \[  A^2_{(2)}(\wt{J_b})\ \xrightarrow{(s_b)_\ast}\ A^2_{(2)}\bigl((S_b)^{[2]}\bigr)\ \xrightarrow{(s_b)^\ast}\ A^2_{(2)}(\wt{J}_b)\ \]
         is the identity.
      Let $\Psi_b\in A^4\bigl( (S_b)^{[2]}\times (S_b)^2\bigr)$ be the correspondence such that $(\Psi_b)^\ast(\Psi_b)_\ast=\ide$ on $A^2_{(2)}\bigl((S_b)^{[2]}\bigr)$. This $\Psi_b$ is obviously the restriction of a relative correspondence $\Psi$ (cf. for instance \cite[Proof of Corollary 3.4]{mmj}).     
       The argument of \cite[Proposition 2.15]{mmj} gives that $A^2_{(2)}\bigl( (S_b)^2\bigr)$ factors (via universally defined correspondences) over $ A^2(S_b)\oplus  A^2(S_b)$. Composing with $\Psi\circ \Gamma_s$ and its transpose, we obtain the required relative correspondences.     
       \end{proof}
       
       Let us now return to the relative correspondence $\Gamma-\gamma\in A^6(\wt{\JJ}\times_B \wt{\JJ})$ constructed above. We define the compositions
       \[ \Gamma_i:=   (\Gamma-\gamma)\circ \Theta_i\ \ \ \in\ A^6(\Ss\times_B \wt{\JJ})\ \ \ (i=1,2)\ .\]
       In view of (\ref{homtriv}), these correspondences are fibrewise homologically trivial:
       \[ (\Gamma_i)\vert_b=0\ \ \ \hbox{in}\ H^{12}(S_b\times\wt{J_b})\ \ \ \forall b\in B\ \ \ (i=1,2)\ .\]
       Applying proposition \ref{gfc4}, it follows that they are fibrewise rationally trivial:
         \[ (\Gamma_i)\vert_b=0\ \ \ \hbox{in}\ A^{6}(S_b\times\wt{J_b})\ \ \ \forall b\in B\ \ \ (i=1,2)\ .\]
     But then a fortiori 
       \[   (\Gamma_i)\vert_b\circ (\Xi_i)\vert_b = (\Gamma-\gamma)\vert_b\circ (\Theta_i)\vert_b\circ (\Xi_i)\vert_b=0\ \ \ \hbox{in}\ A^{6}(\wt{J_b}\times\wt{J_b})\ \ \ \forall b\in B\ \ \ (i=1,2)\ .\]     
       Taking the sum, this implies the fibrewise vanishing
       \[ (\Gamma-\gamma)\vert_b   \circ (\Theta_1\circ\Xi_1+ \Theta_2\circ\Xi_2)\vert_b=0   \ \ \ \hbox{in}\ A^{6}(\wt{J_b}\times\wt{J_b})\ \ \ \forall b\in B\ .\]
       In view of Lemma \ref{ll}, we find that
       \[ \bigl(  (\Gamma-\gamma)\vert_b\bigr){}_\ast=0\ \colon\ \ A^2_{(2)}(\wt{J_b})\ \to\ A^4(\wt{J_b})\ \ \ \forall b\in B\ .\]
       But the correspondence $\gamma\vert_b$ does not act on $A^2(\wt{J_b})$ for dimension reasons, and so
       \[  (\Gamma\vert_b){}_\ast=0\ \colon\ \ A^2_{(2)}(\wt{J_b})\ \to\ A^4(\wt{J_b})\ \ \ \forall b\in B\ .\]
      Since (by construction) $\Gamma\vert_b=\Gamma_b$ acts on $A^2_{(2)}(\wt{J_b})$ as multiplication by $A_b$, this proves claim \ref{c2}.      
     \end{proof}

  \begin{remark} The fourfold $J_0$, being birational to $S^{[2]}$, has an MCK decomposition (theorem \ref{square}).
  In proving theorem \ref{main1}, it would be more natural to use this MCK decomposition of $J_0$, rather than the one of $\wt{J_0}$. However, when trying to do this one runs into the following problem: it is not clear whether the MCK decomposition of $J_0$ is universal (in the sense of lemma \ref{umck}); for this one would need to know that the correspondence $Z$ constructed in \cite{Rie} is universally defined.
  (On a related note, it is not clear whether the map $(r_b)^\ast\colon A^\ast(J_b)\to A^\ast(\wt{J_b})$ respects the bigrading coming from the two MCK decompositions, i.e. I have not been able to prove that $r_b$ is ``of pure grade 0'' in the sense of \cite{SV2}.)
  \end{remark}

 \section{Examples in dimension 6}
 
 \begin{theorem}[Mukai \cite{Mu}]\label{mu} Let $S\subset\PP^3$ be a quartic K3 surface, and assume that every element in $\vert \OO_S(1)\vert$ is irreducible. Let $\pi\colon J_1\to \vert \OO_S(1)\vert\cong\PP^3$ be the component of the compactified Picard scheme that parametrizes torsion free degree $1$ line bundles $\xi$ on curves $C\in \vert \OO_S(1)\vert$. The subset $P\subset J_1$ parametrizing line bundles $\xi$ such that $H^0(C,\xi)\not=0$ has the structure of a $\PP^2$-bundle over $S$. The flop of $P\subset J_1$ is isomorphic to the moduli space $M_v(S)$, where $v$ is the Mukai vector $v=(3,\OO_S(-1),0)$.
 \end{theorem}
 
 \begin{proof} This is \cite[Example 0.8]{Mu}.
 
 \end{proof}

 \begin{theorem}\label{main2} Let $J_1$ and $P$ be as in theorem \ref{mu}. Let $h_1\in A^1(J_1)$ be the polarization class, and let $h_2:=\pi^\ast(d)\in A^1(J_1)$ where $d\subset\PP^3$ is a hyperplane class.
%let $A\in A^3(J_1)$ be a general fibre of the fibration $J_1\to\PP^3$. 
Let $R^\ast(J_1)$ be the $\QQ$-subalgebra 
   \[ R^\ast(J_1):=  \langle  h_1, h_2, P, c_j(T_{J_1})\rangle\ \ \ \ \subset\ A^\ast(J_1)\ .\]    The cycle class map induces an injection
   \[ R^\ast(J_1)\ \hookrightarrow\ H^\ast(J_1,\QQ)\ .\]
 \end{theorem}
 
 \begin{proof} Let $\JJ\to B_{}$ be the universal family of sixfolds $J_1$ as in theorem \ref{mu} (here $B$ is some open in the parameter space $B_{g3}$ of notation \ref{not}). We will prove that $\JJ\to B_{}$ has the Franchetta property. Since the classes defining the subring $R^\ast(J_1)$ are universally defined (i.e. they are restrictions of classes in $A^\ast(\JJ)$), this settles the theorem.
 
 We claim that there exist morphisms of $B$-schemes
   \[  \JJ\ \xleftarrow{r}\ \wt{\JJ}\ \xrightarrow{s}\ \MM\ ,\]
   where $\MM\to B$ is the universal moduli space with Mukai vector $v=(3,\OO_S(-1),0)$, and $r\colon\wt{\JJ}\to\JJ$ is the blow-up of $\PPP$ (the relative version of $P$) and $s\colon\wt{\JJ}\to \MM$ is the blow-up of $\PPP^\prime$ (the relative version of the dual $\PP^2$-bundle $P^\prime\subset M_v$). To ascertain that $\MM$ and $s$ exist as claimed, one may reason as follows: $\PPP\subset\JJ$ can obviously be defined and has the structure of a $\PP^2$-bundle over $\Ss$.
  Let $r\colon\wt{\JJ}\to \JJ$ be the blow-up with center $\PPP$, and let $\EE\subset\wt{\JJ}$ denote the exceptional divisor of $r$.
  This $\EE$ maps to $\PPP^\prime$, which is the dual $\PP^2$-bundle over $Ss$.
   The Nakano--Fujiki criterion for the existence of a blow-down \cite{FN}, as used by Mukai \cite[Proof of Theorem 0.7]{Mu} needs that the normal bundle of
   $\EE\subset\wt{\JJ}$ restricts to the tautological bundle of the fibres $F_p$ of $\EE\to\PPP^\prime$. Since $\NNN_{\EE/\wt{\JJ}}\vert_{F_p}=
   \NNN_{E_b/\wt{\JJ}_b}\vert_{F_p}$, and the criterion is satisfied fibrewise, this is OK. That is, thanks to Nakano--Fujiki we conclude that there exists a blow-down $s\colon\wt{\JJ}\to\MM$ with $\MM$ smooth and $s(\EE)=\PPP^\prime$, as claimed.
   
   To prove the Franchetta property for $\JJ$, it suffices to prove the Franchetta property for $\wt{\JJ}$.
   On the other hand, the morphism $s$ is the blow-up with center $\PPP^\prime$, and $\PPP^\prime$ is a $\PP^2$-bundle over $\Ss$. The formulae for Chow groups of blow-ups and projective bundles give a commutative diagram
   \[
   \begin{array}[c]{ccc}
      A^i(\wt{\JJ}) & \xrightarrow{}& A^i(\MM)\oplus A^{i-1}(\Ss)\oplus  A^{i-2}(\Ss)\oplus     A^{i-3}(\Ss)\\
       &&\\
       \downarrow&&\downarrow\\
       &&\\
      A^i(\wt{J}_b) & \xrightarrow{\cong}& A^i(M_b)\oplus A^{i-1}(S_b)\oplus  A^{i-2}(S_b)\oplus     A^{i-3}(S_b)\ .\\   
      \end{array}\]
 (Here, we write $\wt{J}_b, M_b, S_b$ for the fibre over $b\in B$ of the family $\wt{\JJ}$ resp. $\MM$ resp. $\Ss$.)
 Since we already know the Franchetta property holds for $\Ss\to B$, the Franchetta property for $\wt{\JJ}\to B$ follows from that for $\MM\to B$. Thus, the following result settles the proof of theorem \ref{main2}:
 
 \begin{proposition}\label{p} The family $\MM\to B$ has the Franchetta property.
 \end{proposition}
  
  To prove proposition \ref{p}, we use B\"ulles' result (Theorem \ref{bul}) to reduce to the family $\Ss^{3/B}:=\Ss\times_B \Ss\times_B \Ss$.
  That is, Theorem \ref{bul} tells us that for every $b\in B$ there exist correspondences
    \[ \Gamma_1^b,\ldots,\Gamma_r^b\ \  \in A^\ast(M_b\times (S_b)^{k_j})\ ,\ \ \    \Psi_1^b,\ldots,\Psi_r^b\ \  \in A^\ast( (S_b)^{k_j}   \times M_b)\ \]
    with the property that
    \[ \Delta_{M_b}=\sum_{j=1}^r \Psi_j^b\circ \Gamma_j^b\ \ \ \hbox{in}\ A^6(M_b\times M_b)\ .\]
    
 Using a Hilbert schemes argument as in \cite[Proposition 3.7]{V0}, these fibrewise data can be spread out over the family, i.e. there exist relative correspondences
      \[ \Gamma_1,\ldots,\Gamma_r\ \  \in A^\ast(\MM\times_B \Ss^{k_j/B})\ ,\ \ \    \Psi_1,\ldots,\Psi_r\ \  \in A^\ast( \Ss^{k_j/B}   \times_B \MM)\ \]
    with the property that
    \begin{equation}\label{prop} \Delta_{M_b}=\sum_{j=1}^r (\Psi_j\circ \Gamma_j)\vert_{M_b\times M_b}\ \ \ \hbox{in}\ A^6(M_b\times M_b)\ \ \ \forall b\in B\ .\end{equation}    
(Alternatively, instead of invoking a Hilbert schemes argument, one may observe that the cycles in \cite{Bul} are universal expressions in the Chern classes of a quasi-universal object, and thus naturally can be constructed in families. This is the same argument as \cite[Proof of Theorem 3.1]{FLV2}.)
    
  Now, given a cycle $\Gamma\in A^\ast(\MM)$ which is homologically trivial on the very general fibre, the element
  \[ (\Gamma_1\circ \Gamma,\ldots,\Gamma_r\circ\Gamma)\ \ \in\ A^\ast(\Ss^{k_1/B})\oplus \cdots\oplus   A^\ast(\Ss^{k_r/B})     \]
  will also be homologically trivial on the very general fibre. The families $\Ss^{k_j/B} $ have the Franchetta property by Theorem \ref{gfc} (note that $k_j\le 3$), and so it follows that
    \[ (\Gamma_1\circ \Gamma,\ldots,\Gamma_r\circ\Gamma)\vert_b=(0,\ldots,0)\ \ \hbox{in}\ \ A^\ast(S^{k_1})\oplus \cdots\oplus   
    A^\ast(S^{k_r})  \ \ ,  \]
    for $b\in B$ very general. But then, in view of (\ref{prop}), it follows that also
    \[ \Gamma\vert_b= \sum_{j=1}^r (\Psi_j\circ \Gamma_j\circ \Gamma)\vert_b =0\ \ \ \hbox{in}\ A^\ast(\MM)\ ,\]
    for $b\in B$ very general, i.e. $\MM\to B$ has the Franchetta property.
    
    This proves the proposition, and hence the theorem.
         \end{proof}

  \begin{remark} A general fibre $A$ of the Lagrangian fibration $J_1\to \PP^3$ has $A=h_2^3$ in $A^3(J_1)$ (because a point $p\in\PP^3$ has $p=d^3$ in $A^3(\PP^3)$, with $d$ the hyperplane class). As such, $A$ is in the subalgebra $R^\ast(J_1)$ of Theorem \ref{main2}.
   \end{remark}

 \begin{remark} It would be interesting to prove something like Theorem \ref{main1} for the sixfolds $J_1$ of theorem \ref{main2}. To do this, it would suffice to have an MCK decomposition for $J_1$ that is universal, and with the property that $A^2_{(2)}(J_1)$ comes from $A^2(S)$.
 \end{remark}

\vskip1cm
\begin{nonumberingt} Thanks to Yoyo of kuchibox.fr for daily delivery and great service.

\end{nonumberingt}

\end{document}